\documentclass[reqno,12pt,twoside,a4paper]{amsart} 

\usepackage{bbm,enumerate,amsmath,amstext,mathrsfs,a4}
\usepackage[english]{babel}
\usepackage[latin1]{inputenc}
\usepackage{graphicx}

\newcommand{\R}{\mathbbm{R}}

\newcommand{\Z}{\mathbbm{Z}}
\newcommand{\C}{\mathbbm{C}}
\newcommand{\CP}{\C P}
\newcommand{\RP}{\R P}

\newcommand{\Th}{\mathbbm{T}\mathrm{h}}
\newcommand{\thh}{\mathrm{Th}}
\newcommand{\Diff}{\mathrm{Diff}}
\newcommand{\F}{\mathbbm{F}}

\newcommand{\Tor}{\mathrm{Tor}}
\newcommand{\Cotor}{\mathrm{Cotor}}
\newcommand{\modd}{\backslash\!\!\backslash}

\newcommand{\Ker}{\mathrm{Ker}}
\newcommand{\Cok}{\mathrm{Coker}}
\newcommand{\IM}{\mathrm{Im}}
\newcommand{\moddd}{/\!\!/}

\newcommand{\trf}{\partial}
\renewcommand{\phi}{\varphi}
\renewcommand{\epsilon}{\varepsilon}

\newcommand{\Sq}{\mathrm{Sq}}

\renewcommand{\AA}{\mathscr{A}}

\newcommand{\tensor}{\otimes}
\newcommand{\isom}{\cong}

\newcommand{\colim}{\operatorname{colim}}

\newcommand{\SO}{\mathrm{SO}}
\newcommand{\GL}{\mathrm{Gl}}
\newcommand{\OO}{\mathrm{O}}
\newcommand{\Spin}{\mathrm{Spin}}

\newcommand{\GV}{\mathrm{G}\mathcal{V}}
\newcommand{\GW}{\mathrm{G}\mathcal{W}}
\newcommand{\W}{\mathcal{W}}
\newcommand{\V}{\mathcal{V}}
\newcommand{\U}{\mathcal{U}}
\newcommand{\KK}{\mathscr{K}}
\newcommand{\Sad}{\mathrm{Sad}}
\newcommand{\Wloc}{\mathcal{W}_\mathrm{loc}}
\newcommand{\WlocT}{\mathcal{W}_{\mathrm{loc},T}}

\newcommand{\WT}{\mathcal{W}_T}

\newcommand{\LL}{\mathcal{L}}
\newcommand{\LT}{\mathcal{L}_T}
\newcommand{\Lloc}{\mathcal{L}_\mathrm{loc}}
\newcommand{\LlocT}{\mathcal{L}_{\mathrm{loc},T}}
\newcommand{\Erg}{E^\mathrm{rg}}
\newcommand{\frg}{{f^\mathrm{rg}}}
\newcommand{\hocolim}{\operatorname{hocolim}}
\newcommand{\hofib}{\operatorname{hofib}}
\newcommand{\Bun}{\mathrm{Bun}}
\newcommand{\CM}{\mathcal{C}_M}
\newcommand{\CMt}{\mathcal{C}_{\theta,M}}
\newcommand{\CMop}{\mathcal{C}_M^{\mathrm{op}}}
\newcommand{\CMtop}{\mathcal{C}_{\theta,M}^{\mathrm{op}}}
\newcommand{\MM}{\mathscr{M}}
\newcommand{\spectrumfont}[1]{\mathbf{#1}}

\usepackage{amsthm}
\theoremstyle{plain}
\newtheorem{thm}{Theorem}[section]
\newtheorem{prop}[thm]{Proposition}
\newtheorem{lem}[thm]{Lemma}

\newtheorem{cor}[thm]{Corollary}

\theoremstyle{definition}

\newtheorem{defn}[thm]{Definition}

\newtheorem{example}[thm]{Example}

\theoremstyle{remark}

\newtheorem{rem}[thm]{Remark}

\numberwithin{equation}{section}

\usepackage[ps,all,dvips]{xy}
\CompileMatrices
\SelectTips{cm}{12}

\begin{document}

\title{Mod $2$ homology of the stable spin mapping class group}
\author{Søren Galatius}
\address{Stanford University, Stanford, USA}
\begin{abstract}
  We compute the mod 2 homology of spin mapping class groups in the
  stable range.  In earlier work \cite{G} we computed the stable mod
  $p$ homology of the oriented mapping class group, and the methods
  and results here are very similar.  The forgetful map from the spin
  mapping class group to the oriented mapping class groups induces a
  homology isomorphism for odd $p$ but for $p=2$ it is far from being
  an isomorphism.  We include a general discussion of tangential
  structures on 2-manifolds and their mapping class groups and then
  specialise to spin structures.  As in \cite{MW}, on which \cite{G}
  is based, the stable homology is the homology of the zero space of a
  certain Thom spectrum.
\end{abstract}
\email{galatius@imf.au.dk}

\maketitle\thispagestyle{empty}

\section{Introduction and statement of results}
\label{sec:introduction}

The main result of this paper is a calculation of the mod 2 homology
of the ``spin mapping class groups'' in a stable range, in the spirit
of \cite{G}, which rests heavily on \cite{MW}.  The paper consists of
two parts.  In the second part we adapt the proof in \cite{MW} to the
case of surfaces with a spin structure\footnote{
  The version of \cite{MW} on which this manuscript is based, only
  treats \emph{oriented} surfaces.  A newer version of \cite{MW} that
  treats surfaces with a more general tangential structure, similar to
  the ``$\theta$-structures'' considered in the second part of this
  paper, has since become available.}.  The result is that these
groups have the same homology, in a stable range, as the infinite loop
space $\Omega^\infty\Th(-U_{\Spin(2)})$ where $U_{\Spin(2)}= E\Spin(2)
\times_{\Spin(2)}\R^2$ is the canonical $\Spin(2)$-vectorbundle over
$B\Spin(2)$ and $-U_{\Spin(2)}$ is the $-2$-dimensional virtual
inverse.  $\Th(-U_{\Spin(2)})$ is the Thom spectrum with Thom class in
dimension $-2$ (see section~\ref{sec:cofibration-sequence} for a more
precise definition).  In the first part of the paper we calculate the
mod 2 homology of the infinite loop space
$\Omega^\infty\Th(-U_{\Spin(2)})$.  Let us introduce some notation
before giving a more precise description of our results.

Let $\theta\colon U_3 \to B_3$ be a three-dimensional real vector
bundle, and let $P_3 \to B_3$ be the underlying principal
$\GL_3(\R)$-bundle.  For the moment it can be arbitrary, but we shall
later specialise to the case $\theta = \theta_{\Spin}\colon
E\Spin(3)\times_{\Spin(3)} \R^3 \to B\Spin(3)$.  Let $U_2 = P_3
\times_{\GL_2(\R)} \R^2$.  This is a 2-dimensional real vectorbundle
over the space $B_2 = P_3/\GL_2(\R)$.  $B_2$ is a fibre bundle over
$B_3$ with fibre $\GL_3(\R)/\GL_2(\R) \simeq S^2$, so $B_2$ is fibre
homotopy equivalent to the sphere bundle of $U_3$.  We have a
canonical homeomorphism
\begin{equation*}
  \Bun(V,U_2) = \Bun(V\times \R, U_3)
\end{equation*}
where $\Bun$ denotes the space of bundle maps.

\begin{defn}
  Let $F$ be a surface, possibly with boundary, and let $\theta\colon  U_3
  \to B_3$ be as above.  Then the space of $\theta$-structures on $F$
  is the space $\Bun(TF, U_2)$ of bundle maps (these are suppose to be
  standard near the boundary if $F$ has boundary).  This has a left
  action of $\Diff(F)$.  The space of $(F,\theta)$-surfaces is the
  space
  \begin{equation*}
    \MM^\theta(F) := E\Diff(F) \times_{\Diff(F)} \Bun(TF,U_2).
  \end{equation*}
\end{defn}
The space $\MM^\theta(F)$ is the classifying space for pairs
$(\pi,\xi)$ of a fibre bundle $\pi\colon E \to X$ with fibre $F$ and a
bundle map $\xi\colon T^\pi E \to U_2$, where $T^\pi E$ denotes the
fibrewise tangentbundle of $E$.  Notice that $\Bun(TF,U_2)$ may be
empty.  This will be the case e.g.\ if $U_2$ is orientable but $F$ is
not.  Notice also that $\MM^\theta(F)$ may be non-connected.  We
describe its components.

\begin{sloppypar}
The action of $\Diff(F)$ on $\Bun(TF,U_2)$ induces an action of
$\Diff(F)$ on $\pi_0\Bun(TF,U_2)$.  For $\gamma\in\pi_0\Bun(TF,U_2)$
we write $\Diff(F,\gamma) \subseteq \Diff(F)$ for the subgroup that
fixes $\gamma$.  Define
\begin{equation*}
  \MM^\theta(F,\gamma) = E\Diff(F,\gamma)\times_{\Diff(F,\gamma)}
  \Bun_\gamma(TF,U_2)
\end{equation*}
This is a connected space, and in general we have
\begin{equation*}
  \MM^\theta(F) \simeq \coprod_{\gamma} \MM^\theta(F,\gamma),
\end{equation*}
where the disjoint union is over one $\gamma \in \pi_0\Bun(TF,U_2)$
in each $\Diff(F)$-orbit.
\end{sloppypar}

There is a fibration sequence
\begin{equation}\label{eq:22}
  \Bun_\gamma(TF,U_2) \to \MM^\theta(F,\gamma) \to B\Diff(F,\gamma).
\end{equation}
In particular (for genus $\geq 2$), $\MM^\theta(F,\gamma)$ is a
$K(\pi,1)$ if and only if $\Bun_\gamma(TF,U_2)$ is.  In this case we
have $\MM^\theta(F,\gamma) = B\Gamma^\theta(F,\gamma)$ where
$\Gamma^\theta(F,\gamma) = \pi_1\MM^\theta(F,\gamma)$ is what we could
call the mapping class group of $(F,\gamma)$.

The parametrised Pontryagin-Thom construction defines a map
\begin{equation*}
  \alpha\colon \MM^\theta(F,\gamma) \to \Omega^\infty\Th(-U_2)
\end{equation*}
and in favorable cases this will be an isomorphism in $H_n(-;\Z)$ when
$n$ is small compared to the genus of $F$.

The case $\theta = \theta_{\SO}\colon  E\SO(3)\times_{\SO(3)}\R^3 \to
B\SO(3)$ is equivalent to the case considered in \cite{MW}:  An element
$\gamma \in \pi_0\Bun(TF,U_2)$ is an orientation of $F$, and
$\MM^\theta(F,\gamma) \simeq B\Diff(F,\gamma)$ is the classifying
space of the group of orientation preserving diffeomorphisms.
Furthermore $\Omega^\infty\Th(-U_2) = \Omega^\infty\CP^\infty_{-1}$.
The homology of this space is calculated in \cite{G}.

\begin{sloppypar}
  Now specialise to the case $\theta = \theta_\Spin\colon
  E\Spin(3)\times_{\Spin(3)}\R^3 \to B\Spin(3)$.  Then an element
  $\gamma\in\pi_0\Bun(TF,U_2)$ is a ``spin structure'' on $F$, given
  equivalently by a ``quadratic refinement of the intersection form''
  on $H_1(F,\F_2)$ cf \cite{Johnson}.  Any two spin structures on $F$
  differ by an element in $H^1(F,\F_2)$ so there are $4^g$ spin
  structures.  There are only two $\Diff(F)$-orbits, however.  They
  are distinguished by the Arf invariant of the quadratic form.
  Therefore
\begin{equation*}
  \MM^\theta(F) = \MM^\theta(F,\gamma_0)\amalg \MM^\theta(F,\gamma_1)
\end{equation*}
where $\gamma_0$ is an Arf invariant 0 spin structure and $\gamma_1$
is an Arf invariant 1 spin structure.
\end{sloppypar}

The fibration sequence~\eqref{eq:22} specialises to
\begin{equation*}
  \RP^\infty \to \MM^\theta(F,\gamma) \to B\Diff(F,\gamma)
\end{equation*}
and in case $F$ has genus $\geq 2$ these are all $K(\pi,1)$-spaces.
The fundamental groups of $\MM^\theta(F,\gamma)$ and
$B\Diff(F,\gamma)$ could both be called ``spin mapping class groups''.
Both are studied in \cite{Bauer} who uses the notation $G_\gamma(F) =
\pi_1 B\Diff(F,\gamma) = \pi_0 \Diff(F,\gamma)$ and $\Tilde
G_\gamma(F) = \pi_1 \MM^\theta(F,\gamma)$ and attributes the latter to
Gregor Masbaum.  \cite{Harer} and \cite{Bauer} proves homological
stability of these groups: If $F'$ is obtained from $F$ by glueing
along boundaries of $F$, then the natural maps $G_\gamma(F) \to
G_\gamma(F')$ and $\Tilde G_\gamma(F) \to \Tilde G_\gamma(F')$ are
both isomorphisms in $H_k(-;\Z)$ when $g \geq 2k^2+6k-2$, where $g$ is
the genus of $F$.  \cite{Harer} proves this in the case where $F'$ has
boundary, and \cite{Bauer} extends Harer's proof to the case where
$\partial F' = \emptyset$, and proves homological stability for the
groups $\Tilde G_\gamma(F)$.

The parametrised Pontryagin-Thom construction defines a map
\begin{equation}
  \label{eq:19}
  \alpha\colon  \MM^\theta(F,\gamma) \to \Omega^\infty\Th(-U_{\Spin(2)})
\end{equation}
and we prove that (on components) it is an isomorphism in $H_k(-;\Z)$
whenever $g\geq 2k^2 + 6k-2$.  The number $2k^2 + 6k-2$ is the
stability range of $\Tilde G_\gamma(F)$ in \cite{Bauer}, and an
improvement of the stability range would give an improvement of the
isomorphism range of the map~\eqref{eq:19}.

\begin{sloppypar}
It is easily seen (using e.g.\ the fibration sequence \eqref{eq:1}
below) that $\pi_0 \Omega^\infty \Th( -U_{\Spin(2)} ) \cong \Z
\times\Z/2$.  One may verify that the image of the map~\eqref{eq:19}
is in the component given by the genus of $F$ and the Arf invariant of
$\gamma$.  We conclude this introduction by stating the theorems about
the homology of $\Omega^\infty\Th(-U_{\Spin(2)})$, and thus the
homology in a stable range, of $\MM^\theta(F,\gamma)$.
\end{sloppypar}

The starting point of the calculation is a fibration sequence of
infinite loop spaces
\begin{equation}
  \label{eq:1}
  \xymatrix{
    {\Omega^\infty\Th(-U_{\Spin(2)})} \ar[r]^-{\Omega\omega} & {Q(B\Spin(2)_+)}
    \ar[r]^-{\Omega\partial} & {\Omega Q(S(U_{\Spin(2)})_+)}.
  }
\end{equation}
Here $S(U_{\Spin(2)})$ is the sphere bundle of $U_{\Spin(2)}$ and $Q$
denotes the functor $\Omega^\infty\Sigma^\infty$.  If we identify
$B\Spin(2)$ with $\CP^\infty$ then $U_{\Spin(2)} = L\tensor_\C L$,
where $L$ is the canonical complex line bundle, and $S(U_{\Spin(2)}) =
\RP^\infty$.  We give a concrete description of~\eqref{eq:1} in
Section~\ref{sec:cofibration-sequence}.

For brevity we shall write $U$ for $U_{\Spin(2)}$.  In the following,
all Hopf algebras are commutative and cocommutative.  Recall that any
map $f: A \to B$ of such Hopf algebras have a kernel denoted $A \modd
f$ and a cokernel $B\moddd f$ in the category of Hopf algebras.
Homology and cohomology is always with coefficients in $\F_2$.
\begin{thm}
  \label{thm:1}
  The fibration sequence \eqref{eq:1} induces a short exact sequence
  of Hopf algebras
  \begin{equation}\label{eq:2}
    \xymatrix{
      {H_*(\Omega^\infty\Th(-U)) \modd\Omega\omega_*}
      \ar@{^{(}->}[r]& {H_*(\Omega^\infty\Th(-U)) } \ar@{->>}[r]^-{\Omega\omega_*}&
      {H_*(Q(B\Spin(2)_+))\modd \Omega\trf_*}
    }
  \end{equation}
  and dually
  \begin{equation}\label{eq:3}
    \xymatrix{
      {H^*(Q_0(B\Spin(2)_+))\moddd \Omega\trf_*} 
      \ar@{^{(}->}[r]^-{\Omega\omega^*}& {H^*(\Omega_0^\infty\Th(-U)) } \ar@{->>}[r]&
      {H^*(\Omega_0^\infty\Th(-U)) \moddd\Omega\omega^*.}
    }
  \end{equation}
\end{thm}

\begin{sloppypar}
It remains to determine the Hopf algebras $H_*(Q(B\Spin(2)_+))\modd
\Omega\trf_*$ and $H_*(\Omega^\infty\Th(-U)) \modd\Omega\omega_*$.
The next theorem determines the Hopf algebra
$H_*(Q(B\Spin(2)_+))\modd\Omega\trf_*$.  We also produce an explicit
splitting of the sequence~\eqref{eq:2}, although the splitting is only
as algebras, not as Hopf algebras.
\end{sloppypar}


The action of $\Spin(3) = SU(2)$ on $S^2$ gives an $S^2$-bundle
$E\Spin(3)\times_{\Spin(3)} S^2 \to B\Spin(3)$.  The vertical tangent
bundle $E\Spin(3)\times_{\Spin(3)} TS^2$ has a canonical
spin-structure, and the classifying map $E\Spin(3) \times_{\Spin(3)}
S^2 \to B\Spin(2)$ is a homotopy equivalence.  Consequently we get a
map $B\Spin(3) \to \MM^\theta(S^2)$.  The composition
\begin{equation*}
  B\Spin(3) \to \MM^\theta(S^2) \stackrel{\alpha}{\to} \Omega^\infty \Th(-U) \to
  Q(B\Spin(2)_+)
\end{equation*}
is the Becker-Gottlieb transfer for the fibration sequence $S^2 \to
B\Spin(2) \to B\Spin(3)$.

Before stating the next theorem, let us recall that for a Hopf algebra
$A$ over $\F_2$ there is a Frobenius map $\xi\colon  A \to A$ given by $\xi
x = x^2$ which is a morphism of Hopf algebras.  Write $a_i\in
H_{2i}(B\Spin(2))$ and $b_i \in H_{4i}(B\Spin(3))$ for the generators,
$i\geq 0$.  Recall that $H_*(Q(B\Spin(2)_+))$ is the free commutative
algebra on the set $\mathbf{T}_2$ of generators given by
\begin{equation*}
  \mathbf{T}_2 = \{ Q^I a_i \mid \text{$i\geq 0$, $I$ admissible,
    $e(I) > 2i$}\},
\end{equation*}
where $Q^I$ are the iterated Dyer-Lashof operations (see \cite{CLM}
for definitions and proofs).  Similarly $H_*(Q(B\Spin(3)_+))$ is the
free commutative algebra on the set of generators given by
\begin{equation*}
  \mathbf{T}_3 = \{Q^I b_i \mid \text{$i\geq 0$, $I$ admissible,
    $e(I) > 4i$}\}.
\end{equation*}

\begin{thm}\ 
  \label{thm:2}
  \begin{enumerate}[(i)]
  \item We have $H_*(Q(B\Spin(2)_+))\modd\trf_* = \xi
    H_*(Q(B\Spin(2)_+))$.  Both the algebra
    $H_*(Q(B\Spin(2)_+))\modd\Omega\trf_*$ and
    the dual algebra $H^*(Q_0(B\Spin(2)_+))\moddd\Omega\trf^*$ are free
    commutative.
  \item   The composition
    \begin{equation*}
      H_*(Q(B\Spin(3)_+)) \to H_*(\Omega^\infty\Th(-U)) \to
      H_*(Q(B\Spin(2)_+)) \modd \Omega\trf_*
    \end{equation*}
    is surjective.  It maps $b_i$ to $a_i^2$ and more generally it
    maps $Q^{2I}b_i$ to $(Q^I a_i)^2$.
  \end{enumerate}
\end{thm}

It remains to describe the Hopf algebra $H_*(\Omega^\infty\Th(-U))
\modd\Omega\omega_*$ in Theorem~\ref{thm:1}.  This is done by first
describing the (co-)homology of $\Omega Q(\RP^\infty_+)$ and $\Omega^2
Q(\RP^\infty_+)$.

To state the results about $\Omega Q(\RP^\infty_+)$ and $\Omega^2
Q(\RP^\infty_+)$, let us recall a certain functor from \cite{MM}.  It
is called $V$ in \cite[definition 6.2]{MM}, but we shall call it $A$.
We shall as usual let $P(-)$ and denote the vectorspace of primitive
elements and $Q(-)$ denote the vectorspace of indecomposable elements
in a Hopf algebra.
\begin{defn}[{\cite{MM}}]
  Let $V$ be a graded vectorspace and $\xi\colon  V \to V$ a linear map such
  that $\xi V_n \subseteq V_{2n}$.  Let $SV$ denote the free
  commutative (i.e.\ polynomial) algebra generated by $V$, and let
  $I\subseteq SV$ be the ideal generated by the elements $x^2 - \xi
  x$, $x\in V$.  Let $AV = A(V,\xi) = SV/I$.
\end{defn}
The functor $A$ satisfies $A(V\oplus V') = AV \tensor AV'$ and
therefore the diagonal $V \to V\oplus V$ induces a comultiplication on
$AV$ making it a Hopf algebra.  The vectorspace of primitive elements
is $V$ itself, $PAV = V$.

\begin{thm}\ 
  \label{thm:3}
  \begin{enumerate}[(i)]
  \item The suspension
    \begin{equation*}
      \sigma^*\colon  QH^*(Q_0\RP^\infty_+) \to PH^*(\Omega Q\RP^\infty_+)
    \end{equation*}
    is an isomorphism (of degree $-1$).
  \item The suspension $\sigma^*$ above induces an isomorphism
    \begin{equation*}
      A( s^{-1} QH^*(Q_0\RP^\infty_+), s^{-1} \Sq_1) \cong
      H^*(\Omega Q\RP^\infty_+).
    \end{equation*}
    Here $s^{-1}$ denotes desuspension of graded vector spaces and
    $\Sq_1$ is the Steenrod operation given by $\Sq_1(x) =
    \Sq^{k-1}(x)$ if $\deg(x) = k$.
  \item The Hopf algebra $H^*(\Omega_0 Q(\RP^\infty_+))$ is
    primitively generated and polynomial.
  \item The suspension induces an isomorphism
    \begin{equation*}
      \sigma^*\colon  \Cok(\Sq_1) \to QH^*(\Omega Q\RP^\infty_+).
    \end{equation*}
  \end{enumerate}
\end{thm}
\begin{thm}\ 
  \label{thm:4}
  \begin{enumerate}[(i)]
  \item The suspension
    \begin{equation*}
      \sigma^*\colon  QH^*(\Omega_0Q\RP^\infty_+) \to PH^*(\Omega^2 Q\RP^\infty_+)
    \end{equation*}
    is an isomorphism (of degree $-1$).
  \item The suspension $\sigma^*$ above induces an isomorphism
    \begin{equation*}
      A( s^{-2} \Cok(\Sq_1), s^{-2} \Sq_2) \cong
      H^*(\Omega^2 Q\RP^\infty_+).
    \end{equation*}
    Here $\Sq_2\colon  \Cok(\Sq_1) \to \Cok(\Sq_1)$ is the Steenrod
    operation given by $\Sq_2(x) = \Sq^{k-2}(x)$ if $\deg(x) = k$.
  \item The Hopf algebra $H^*(\Omega_0 Q(\RP^\infty_+))$ is
    primitively generated but not polynomial.
  \item The suspension induces an isomorphism
    \begin{equation*}
      \sigma^* \circ \sigma^*\colon  \Cok(\Sq_2) \to QH^*(\Omega^2 Q\RP^\infty_+)      
    \end{equation*}
    of degree $-2$.
  \end{enumerate}
\end{thm}

\begin{sloppypar}
Using this description of $H^*(\Omega^2 Q\RP^\infty_+)$ we describe
the Hopf algebra $H_*(\Omega^\infty\Th(-U)) \modd \omega_*$ and its
dual $H^*(\Omega^\infty\Th(-U)) \moddd \omega^*$.
\begin{thm}\ 
  \label{thm:5}
  \begin{enumerate}[(i)]
  \item The Hopf algebra $H_*(\Omega^\infty \Th(-U)) \modd
    \Omega\omega_*$ is precisely the image of $H_*(\Omega^2
    Q\RP^\infty) \to H_*(\Omega^\infty \Th(-U))$.
  \item $H^*(\Omega_0^\infty \Th(-U)) \moddd \Omega\omega^*$ injects
    into $H^*(\Omega^2_0 Q(\RP^\infty_+))$ and is primitively
    generated. 
  \item Under the isomorphism in Theorem~\ref{thm:4}.(ii),
    $H^*(\Omega^\infty_0\Th(-U))\moddd\Omega\omega^*$ is precisely the
    subalgebra generated by the double suspension of the sub
    vectorspace
    \begin{equation*}
      \Ker\left( Q\trf^*\colon  QH^*(Q\RP^\infty_+) \to QH^*(Q\Sigma(
      B\Spin(2)_+)) \right)
    \end{equation*}
    of $QH^*(Q\RP^\infty_+)$.
  \end{enumerate}
\end{thm}
\end{sloppypar}

Finally we can combine the above to conclude the following corollary.
\begin{cor}
  \label{cor:1.8}
  The infinite loop map
  \begin{equation*}
    \Omega^2 Q(\RP^\infty_+) \times Q(B\Spin(3)_+) \to \Omega^\infty\Th(-U),
  \end{equation*}
  which on the first factor is the map $\Omega^2 Q(\RP^\infty_+) \to
  \Omega^\infty\Th(-U)$ induced by~\eqref{eq:1} and which on the
  second factor is the map $Q(B\Spin(3)_+) \to \Omega^\infty\Th(-U)$
  from Theorem~\ref{thm:2}, induces an injection
  \begin{equation*}
    H^*(\Omega^\infty\Th(-U)) \to H^*(\Omega^2 Q(\RP^\infty_+))
    \tensor H^*(Q(B\Spin(3)_+)).
  \end{equation*}
\end{cor}

\subsection{Acknowledgements} This paper is part of my thesis at the
University of Aarhus.  It is a great pleasure to thank my thesis
advisor Ib Madsen for his help and encouragement during my years as a
graduate student.  I also thank M.\ B\"okstedt, J.\ Tornehave and N.\
Wahl for many useful conversations and Lars Madsen for excellent
technical assistance.

\section{A cofibration sequence}
\label{sec:cofibration-sequence}

Let us first describe a concrete model for the maps of spectra
underlying the fibration sequence~\eqref{eq:1}.

Let $q\colon  \CP^n \to \CP^n$ denote the map $q([z_0: \dots : z_n]) =
[z_0^2: \dots : z_n^2]$.  Let $L_n$ denote the canonical complex line
bundle over $\CP^n$ and $L_n^\perp$ its orthogonal complement.  There
is a bundle map
\begin{equation*}
  \xymatrix{
    {L_n\tensor L_n} \ar[r]^{\hat q} \ar[d] & {L_n}\ar[d]\\
    {\CP^n}\ar[r]^q & {\CP^n}\\
  }
\end{equation*}
where $\hat q \colon  (z_0, \dots, z_n) \tensor (w_0, \dots, w_n) \mapsto
(z_0w_0, \dots, z_nw_n)$.  Thus $\hat q$ identifies $L_n \tensor L_n$
with $q^* L_n$.  We shall write $L_n^2 = q^* L_n$ and ${L_n^2}^\perp =
q^* L_n^\perp$.

There is an obvious bundle map
\begin{equation*}
  \xymatrix{
    {L_{n-1}^\perp \times \C} \ar[r]\ar[d] & {L_n^\perp} \ar[d] \\
    {\CP^{n-1}} \ar[r] & {\CP^n}
  }
\end{equation*}
and an induced bundle map
\begin{equation*}
  \xymatrix{
    {{L_{n-1}^2}^\perp \times \C} \ar[r] \ar[d] & {{L_n^2}^\perp} \ar[d] \\
    {\CP^{n-1}} \ar[r] & {\CP^n}
  }
\end{equation*}
These gives maps of Thom spaces $\thh(L_{n-1}^\perp) \wedge S^2 \to
\thh(L_n^\perp)$ and $\thh({L_{n-1}^2}^\perp) \wedge S^2 \to
\thh({L_n^2}^\perp)$.  Therefore we get spectra $\Th(-L)$ and
$\Th(-L^2)$ with $(2n+2)$-nd space $\thh(L_n^\perp)$ and
$\thh({L_n^2}^\perp)$, respectively.  The associated infinite loop
spaces are
\begin{gather}
  \label{eq:4}
  \begin{split}
    \Omega^\infty\Th(-L) &= \colim \Omega^{2n+2} \Th(L_n^\perp)\\
    \text{and}\quad\Omega^\infty\Th(-L^2) &= \colim \Omega^{2n+2}
    \Th({L^2_n}^\perp)
  \end{split}
\end{gather}

The bundle $L \to \CP^\infty$ above is isomorphic to $U_{\SO(2)} =
E\SO(2) \times_{\SO(2)} \R^2 \to B\SO(2)$ and $L^2 \to \CP^\infty$ is
isomorphic to $U_{\Spin(2)} = E\Spin(2) \times_{\Spin(2)} \R^2\to
B\Spin(2)$.  The map $q$ above is induced from the double cover
$\Spin(2) \to \SO(2)$.  Therefore we shall write
$\Omega^\infty\Th(-U_{\SO(2)})$ and $\Omega^\infty\Th(-U_{\Spin(2)})$
for the spaces of~\eqref{eq:4}.

For a vector bundle $\xi \to X$, let
$\thh(\xi) = \xi\cup\{\infty\}$ be the one-point compactification of
the total space.

\begin{lem}\label{lem:2.1}
  Let $\xi$ and $\eta$ be vector bundles over $X$.  Then there is a
  cofibration sequence
  \begin{equation*}
    \xymatrix{
      {\thh(\xi)} \ar[r]^-{z} & {\thh(\xi\oplus\eta)}\ar[r]^-{\trf} &
      {\thh(\R\oplus\xi_{|S(\eta)})}
    }
  \end{equation*}
  where $z$ is induced from the zero section of $\eta$ and
  $\xi_{|S(\eta)}$ denotes pullback of $\xi$ to the sphere bundle of
  $\eta$.  If $\xi\oplus \eta = \R^n \times X$, then $\partial$ is the
  parametrised Pontryagin-Thom construction of the sphere bundle
  $S(\eta) \to X$.
\end{lem}
\begin{proof}
  The normal bundle of the embedding $S(\eta) \to \eta$ is $\R\times
  S(\eta)$.  This embeds via ``polar coordinates'' onto $\eta - X$.
  Therefore the normal bundle of the composition $S(\eta) \to \eta \to
  \eta\oplus\xi$ is $\R\oplus\xi_{|S(\eta)}$ and this embeds onto
  $\xi\oplus\eta - \xi \subseteq \xi\oplus\eta$.  This defines a
  homeomorphism
  \begin{equation*}
    \thh(\xi\oplus\eta)/\thh(\eta) = (\xi\oplus\eta -
    \eta)\cup\{\infty\} \cong \thh(\R\oplus\xi_{|S(\eta)})
  \end{equation*}
  If $\xi\oplus\eta = X\times\R^n$, then $\partial$ is exactly the
  Thom-Pontryagin construction applied to the embedding $S(\eta)
  \subseteq X\times\R^n$ over $X$.
\end{proof}
\begin{lem}
  \label{lem:2.2}
  The map $\RP^{2n+1} \to \CP^n \times \C^{n+1}$ given by
  \begin{equation*}
    [x_0: y_0: \dots : x_n: y_n] \mapsto ([z_0: \dots : z_n],(z_0^2 ,
    \dots, z_n^2)),
  \end{equation*}
  where $z_j = x_j + iy_j$, is a homeomorphism onto $S(L_n^2)
  \subseteq \CP^n\times \C^{n+1}$.  Thus $S(L_n^2) \to \CP^n$ is
  identified with the quotient map
  \begin{equation*}
    \RP^{2n+1} = S^{2n+1}/\{\pm 1\} \to S^{2n+1}/S^1 = \CP^n
  \end{equation*}
  \qed
\end{lem}
\begin{cor}\label{cor:2.3}
  There is a cofibration sequence
  \begin{equation}
    \label{eq:5}
    \xymatrix{
      {\thh(\R\oplus {L_n^2}^\perp)} \ar[r]^-{z} &
      {\Sigma^{2n+3}\CP^n_+} \ar[r]^-{\trf} & {\Sigma^{2n+2}
        \RP^{2n+1}_+}
    }
  \end{equation}
\end{cor}
\begin{proof}
  Let $\xi = \R\oplus{L_n^2}^\perp$ and $\eta = {L_n^2}$ in
  Lemma~\ref{lem:2.1}.  Then $\xi \oplus \eta =
  \CP^n\times\C^{n+1}\times\R$ and $\R\oplus\xi_{|S(\eta)} =
  \C\oplus{L_n^2}^\perp_{|S(L_n^2)} =
  L_n^2\oplus{L_n^2}^\perp_{|S(L_n^2)} = S(L_n^2) \times \C^{n+1}$,
  using the canonical trivialisation of ${L_n^2}_{|S(L_n^2)}$.  Now
  lemmas~\ref{lem:2.1} and~\ref{lem:2.2} gives the desired result.
\end{proof}
\begin{cor}
  \label{cor:2.4}
  There is a cofibration sequence of spectra
  \begin{equation*}
    \xymatrix{
      {\Sigma\Th(-L^2)} \ar[r] & {\Sigma^{\infty+1}(\CP^\infty_+)}
      \ar[r] & {\Sigma^{\infty}\RP^\infty_+}
    }
  \end{equation*}
  and associated fibration sequences
  \begin{equation}\label{eq:6}
    \xymatrix{
      {\Omega^\infty\Sigma\Th(-L^2)} \ar[r] ^-\omega & {Q\Sigma(\CP^\infty_+)}
      \ar[r]^-{\trf} & {Q\RP^\infty_+}
    }
  \end{equation}
  and 
  \begin{equation}\label{eq:7}
    \xymatrix{
      {\Omega^\infty\Th(-L^2)} \ar[r] ^-{\Omega\omega} & {Q(\CP^\infty_+)}
      \ar[r]^-{\Omega\trf} & {\Omega Q\RP^\infty_+}
    }
  \end{equation}
\end{cor}
\begin{prop}
  \label{prop:2.5}
  The map
  \begin{equation*}
    \trf\colon  Q\Sigma\CP^\infty_+ \to Q\RP^\infty_+
  \end{equation*}
  is the ``$S^1$-transfer'' denoted $t_1$ in \cite{MMM}
\end{prop}
\begin{proof}
  The map $t_1$ in \cite{MMM} is exactly the pretransfer of the
  $S^1$-bundle $ES^1 \times_{S^1} (S^1/\{\pm 1\}) \to BS^1$, and this
  is $\trf$.
\end{proof}
\begin{thm}[{\cite{MMM}}]
  \label{thm:2.6}
  Let $\Bar a_r \in H_{2r+1}(\Sigma\CP^\infty_+)$ and $e_r \in
  H_r(\RP^\infty)$ be the generator.  Then
  \begin{equation*}
    \trf_*(\Bar a_r) \equiv e_{2r+1} + Q^{r+1}e_r
  \end{equation*}
  modulo decomposable elements.
\end{thm}
\begin{proof}
  This follows from \cite[Theorem~4.4]{MMM} by ignoring the decomposable terms.
\end{proof}
\begin{cor}
  \label{cor:2.7}
  The map
  \begin{equation*}
    \trf_*\colon  H_*(Q\Sigma\CP^\infty_+) \to H_*(Q\RP^\infty_+)
  \end{equation*}
  is injective.
\end{cor}
\begin{proof}
  This follows from Theorem~\ref{thm:2.6} and the known structure of
  $H_*(Q\Sigma\CP^\infty_+)$ and $H_*(Q\RP^\infty_+)$, cf \cite{CLM}.
\end{proof}

\section{Cohomology of $\Omega Q\RP^\infty_+$ and
  $\Omega^2Q\RP^\infty_+$}

The goal of this section is to prove Theorems~\ref{thm:3} and
\ref{thm:4}.  This is done via the following proposition.

\begin{prop}
  \label{prop:3.1}
  Let $X$ be a simply connected, homotopy commutative, homotopy
  associative $H$-space.  Assume that $H_*(X)$ and $H_*(\Omega
  X)$ are of finite type.  Then  $H^*(X)$ is a polynomial algebra is
  and only if $\xi\colon  PH^*(X) \to PH^*(X)$ is injective.  In this case
  we have
  \begin{enumerate}[(i)]
  \item The suspension
    \begin{equation*}
      \sigma^*\colon  QH^*(X) \to PH^*(\Omega X)
    \end{equation*}
    is an isomorphism (of degree $-1$).
  \item The suspension $\sigma^*$ above induces an isomorphism
    \begin{equation*}
      A[ s^{-1} QH^*(X), s^{-1}\Sq_1] \cong H^*(\Omega X).
    \end{equation*}
    Here $s^{-1}$ denotes desuspension of graded vectorspaces and
    $Sq_1\colon  QH^*(X) \to QH^*(X)$ is the Steenrod operation given by
    $\Sq_1(x) = \Sq^{k-1} (x)$ if $\deg(x) = k$.
  \item The Hopf algebra $H^*(\Omega X)$ is primitively generated.  It
    is polynomial if and only if $\Sq_1\colon  QH^*(X) \to QH^*(X)$ is
    injective.
  \end{enumerate}
\end{prop}
\begin{proof}
  It follows from Borel's structure theorem that $H^*(X)$ is
  polynomial if and only if $\xi\colon  H^*(X)\to X^*(X)$ is injective.  And
  this happens if and only if $\xi\colon  PH^*(X) \to PH^*(X)$ is
  injective.  The proof of the proposition is based on the
  Eilenberg-Moore spectral sequence, see \cite{EM} or the review in \cite{G}.
  The $E_2$-term is $\Tor_{H^*(X)}(\F_2,\F_2)$ and it converges to
  $H^*(\Omega X)$.

  When $H^*(X)$ is a polynomial algebra, the $E_2$-term of the
  spectral sequence is
  \begin{equation*}
    \Tor_{H^*(X)}(\F_2,\F_2) = E[ s^{-1} QH^*(X) ]
  \end{equation*}
  which has generators and primitives concentrated on the line
  $E_2^{-1,*}$.  Therefore it must collapse, because it is a spectral
  sequence of Hopf algebras.  The suspension can be identified with
  the map
  \begin{equation*}
    QH^*(X) \cong \Tor^{-1,*}_{H^*(X)}(\F_2,\F_2) = E_2^{-1,*} \to
    E_2^{-1,*} \subseteq \Tilde H^*(\Omega X)
  \end{equation*}
  and the image is within the vectorspace of primitive elements.
  Therefore $\sigma^*$ is injective because $E^2 = E^\infty$.
  
  The image of $\sigma^*$ generates the algebra $H^*(\Omega X)$
  because it generates the $E^\infty$-term.  In particular we have
  proved that $H^*(\Omega X)$ is primitively generated.

  The suspension $\sigma^*$ commutes with Steenrod operations.  In
  particular we have
  \begin{equation*}
    (\sigma^*(x))^2 = \sigma^*(\Sq_1 x)
  \end{equation*}
  so the image of $\sigma^*$ is closed under the Frobenius map $\xi\colon  x
  \mapsto x^2$.  That $\sigma^*$ is surjective now follows from the
  Milnor-Moore exact sequence,
  \begin{equation*}
    0 \to P\xi H^*(\Omega X) \to PH^*(\Omega X) \to QH^*(\Omega X) \to
    0.
  \end{equation*}
  Namely, if $\sigma^*$ were not surjective, there would be an element
  of minimal degree in $PH^*(\Omega X)$ not in the image of
  $\sigma^*$.  This element would have to map to zero in $QH^*(\Omega
  X)$ because the image of $\sigma^*$ generates.  Hence, by the exact
  sequence, it would have to be a square of some other element.  But
  this contradicts minimality because the image of $\sigma^*$ is
  closed under $\xi$.

  We have proved (i) and the first part of (iii).  Now (ii) follows
  from the fact that $H^*(\Omega X)$ is primitively generated and that
  $\xi\colon  PH^*(\Omega X) \to PH^*(\Omega X)$ corresponds under
  $\sigma^*$ to $\Sq_1$.  Finally, by (ii) we have that $\xi\colon 
  H^*(\Omega X) \to H^*(\Omega X)$ is injective if and only if $\Sq_1\colon 
  QH^*(X) \to QH^*(X)$ is injective.
\end{proof}
\begin{rem}
  Without the assumption on simple connectivity the above proposition
  is generally false.  It does hold in the following very special
  case, however.  Namely, if $\pi_1 X$ is an $\F_2$-vectorspace and
  $X$ splits as $X \simeq \Tilde X \times B\pi_1 X$.  In this case we
  have $PH_1(X) = H_1(X) = \pi_1(X)$ and $QH_0(\Omega X) = \pi_0
  (\Omega X) = \pi_1 (X)$, and for $k\geq 2$ we have $PH_k(X) =
  PH_k(\Tilde X)$ and $QH_{k-1}(\Omega X) = QH_{k-1}(\Omega\Tilde
  X)$.
\end{rem}

Let $e_r \in H_r(\RP^\infty)$ be the generator.  Recall from \cite{CLM}
that $H_*(Q\RP^\infty_+)$ is the free commutative algebra on the set
\begin{equation*}
  \mathbf{T} = \{ Q^I e_r \mid \text{$r\geq 0$, $I$ admissible, $e(I)
  > r$} \}
\end{equation*}

We shall also need a basis for $PH_*(Q\RP^\infty_+)$
\begin{defn}
  \label{defn:3.3}
  Let $p_{2r+1} \in PH_*(Q\RP^\infty_+)$ be the unique primitive class
  with $p_{2r+1} - e_{2r+1}$ decomposable.  For an admissible sequence
  of the form $I = (2s+1, 2I')$ with $e(I) \geq 2i$, let $p_{(I,2i)}$
  be the unique primitive class with $p_{(I,2i)} - Q^I e_{2i}$
  decomposable.  For an admissible sequence $I = (I', 2s+1, 2I'')$
  with $e(I) \geq 2i$, let $p_{(I,2i)} = Q^{I'} p_{(2s+1, 2I'', 2i)}$.
\end{defn}
Thus $p_{(I,i)} \in PH_*(Q\RP^\infty_+)$ is defined for alle $(I,i)$
with $2 \not| (I,i)$.
\begin{lem}
  \label{lem:3.4}
  The set
  \begin{equation*}
    \{ p_{(I,i)} \mid \text{ $i\geq 0$, $I$ admissible, $e(I) \geq i$,
    $2\not| (I,i)$}\}
  \end{equation*}
  is a basis of $PH_*(Q\RP^\infty_+) = PH_*(Q_0\RP^\infty_+)$.
\end{lem}
\begin{proof}
  This is well known.  That $p_{(I,i)}$ spans all of
  $PH_*(Q\RP^\infty)$ follows from the Milnor-Moore exact sequence.
  See \cite{G} for more details
\end{proof}
\begin{defn}
  \label{defn:3.5}
  Define operations $H_*(Q\RP^\infty) \to H_*(Q\RP^\infty)$ by
  \begin{align*}
    \lambda x &= \Sq^k _* x, \quad \deg(x) = 2k,\\
    \lambda' x &= \Sq^k _* x, \quad \deg(x) = 2k+1\\
    \lambda'' x &= \Sq^k _* x, \quad \deg(x) = 2k+2\\
  \end{align*}
  We write $\lambda$, $\lambda'$ and $\lambda''$ for the induced
  operations on $PH_*(Q\RP^\infty_+)$ and $QH_*(Q\RP^\infty_+)$ also.
  These are dual to $\xi = \Sq_0$, $\Sq_1$, and $\Sq_2$ on cohomology,
  respectively.
\end{defn}
\begin{lem}
  \label{lem:3.6}
  In $H_*(\RP^\infty)$ we have
  \begin{align*}
    \lambda e_{2r} &= e_r\\
    \lambda' e_{2r-1} &= r e_r \\
    \lambda'' e_{2r-2} &= \binom r2 e_r\\
  \end{align*}
\end{lem}
\begin{proof}
  This is dual to the formula $\Sq^k w_1^n = \binom nk w_1^{n+k} \in
  H^*(\RP^\infty)$.
\end{proof}
\begin{lem}
  \label{lem:3.7}
  The operations $\lambda$, $\lambda'$ and $\lambda''$ satisfiy the
  relations
  \begin{align}
    \label{eq:8}
    \lambda Q^{2s}x &= Q^s \lambda x \\
    \label{eq:10}
    \lambda' Q^{2s}x &= Q^s \lambda'x \\
    \label{eq:11}
    \lambda' Q^{2s-1}x &= (\deg Q^s \lambda x) Q^s \lambda x \\
    \label{eq:12}
    \lambda'' Q^{2s} x &= Q^s \lambda'' x, \qquad \text{if $\lambda x =
      0$} \\
    \label{eq:13}
    \lambda'' Q^{2s-1} x &= (1+\deg Q^s\lambda' x)Q^s\lambda' x
  \end{align}
\end{lem}
\begin{proof}
  This follows from the Nishida relations (cf \cite{CLM}).
\end{proof}
\begin{prop}\label{prop:3.8}
  $\lambda\colon  QH_*(Q\RP^\infty_+) \to QH_*(Q\RP^\infty_+)$ is
  surjective.
\end{prop}
\begin{proof}
  This is because $\lambda\colon  H_*(\RP^\infty) \to H_*(\RP^\infty)$ is
  surjective.  Explicitly,~\eqref{eq:8} and Lemma~\ref{lem:3.6}
  implies that
  \begin{equation*}
    \lambda Q^{2I}e_{2r} = Q^I e_r
  \end{equation*}
  so the basis $\mathbf{T}$ of $QH_*(Q\RP^\infty_+)$ is hit.
\end{proof}
\begin{prop}\label{prop:3.9}
  $\lambda'\colon  PH_*(Q\RP^\infty_+) \to PH_*(Q\RP^\infty_+)$ is
  surjective.
\end{prop}
\begin{proof}
  Lemma~\ref{lem:3.6} and Lemma~\ref{lem:3.7} imply that
  \begin{equation*}
    \lambda' e_{4r+1} = e_{2r+1}
  \end{equation*}
  and that
  \begin{equation*}
    \lambda' (Q^{4s+1} Q^{4I'} e_{2i}) = Q^{2s+1}Q^{2I'} e_i
  \end{equation*}
  Hence, since $\lambda'$ preserves decomposables
  \begin{equation*}
    \lambda' p_{4r+1} = p_{2r+1}
  \end{equation*}
  and
  \begin{equation*}
    \lambda' p_{(4s+1, 4I', 4i)} = p_{(2s+1, 2I', 2i)}
  \end{equation*}
  and hence
  \begin{equation*}
    \lambda' p_{(2I',4s+1, 4I'', 4i)} = p_{(I',2s+1, 2I'', 2i)}
  \end{equation*}
  Therefore, by Lemma~\ref{lem:3.4}, $\lambda'$ is surjective.
\end{proof}
\begin{prop}
  \label{prop:3.10} $\lambda''\colon  PH_*(Q\RP^\infty_+) \to
  \Ker(\lambda')$ is not surjective.
\end{prop}
\begin{proof}
  The element $p_3 = e_3 + e_1e_2 + e_1^3$ satisfies $\lambda'(p_3) =
  Q^2 e_1 = p_{(1,1)}$ and the element $p_{(2,1)}$ satisfies $\lambda'
  p_{(2,1)} = p_{(1,1)}$.  So $p_{(2,1)} + p_3 \in \Ker(\lambda')$.
  But $PH_4(Q\RP^\infty)$ has basis $\{Q^3 e_1, Q^2 Q^1 e_1\}$ and
  $\lambda'' (Q^3 e_1) = \lambda'' (Q^2 Q^1 e_1) = 0$, so $p_{(2,1)} +
  p_3$ is not hit by $\lambda''$.
\end{proof}

\begin{proof}[Proof of Theorem~\ref{thm:3}]
  This follows from Proposition~\ref{prop:3.1}, using
  Propositions~\ref{prop:3.8} and \ref{prop:3.9}.
\end{proof}
\begin{proof}[Proof of Theorem~\ref{thm:4}]
  This follows from Proposition~\ref{prop:3.1}, using
  Theorem~\ref{thm:3}.

  That $H^*(\Omega^2_0 Q(\RP^\infty_+))$ is not polynomial follows
  from proposition~\ref{prop:3.10}.  Indeed there must be an a
  generator of degree two with square zero.
\end{proof}

\section{The spectral sequence}
\label{sec:spectral-sequence}

The aim of this section is to prove theorems~\ref{thm:1} and
\ref{thm:5}.  The starting point is the fibration~\eqref{eq:1}.  None
of the spaces in the fibration are connected.  In fact we have
\begin{equation*}
  \pi_0 (\Omega Q\RP^\infty_+) = \Z/2\times\Z/2, \quad \pi_0
  Q(B\Spin(2)_+) = \Z, \quad \pi_0 \Omega^\infty\Th(-U) = \Z\times \Z/2
\end{equation*}
and
\begin{equation*}
  \pi_1 Q(B\Spin(2)_+) = \Z/2, \quad \pi_1(\Omega Q(\RP^\infty_+)) =
  \Z/2 \times \Z/2
\end{equation*}

The claim in theorem~\eqref{eq:1} is clearly equivalent to the claim
that the sequence
\begin{equation}
  \label{eq:20}
  \xymatrix{
    {H_*(\Omega^\infty\Th(-U))} \ar[r]^-{\Omega \omega_*} &
    {H_*(Q(B\Spin(2)_+))} \ar[r]^-{\Omega\trf_*} &
    {H_*(\Omega Q(\RP^\infty_+))}
  }
\end{equation}
is short exact (both means that $\Omega\omega_*$ maps onto the kernel
of $\Omega\trf_*$).  This is equivalent to proving that the sequence
\begin{equation}
  \label{eq:21}
  \xymatrix{
    {H_*(\Omega_0^\infty\Th(-U))} \ar[r]^-{\Omega \omega_*} &
    {H_*(Q_0(B\Spin(2)_+))} \ar[r]^-{\Omega\trf_*} &
    {H_*(\Hat\Omega_0 Q(\RP^\infty_+))}
  }
\end{equation}
is short exact.  Here $\Hat\Omega_0 Q(\RP^\infty_+)$ is the double
cover of $\Omega_0 Q(\RP^\infty_+)$ corresponding to the image of
$\Omega\trf$ in $\pi_1$.  This is equivalent because there is a map
from~\eqref{eq:21} to~\eqref{eq:20}, the kernel of which is the
sequence
\begin{equation*}
  \xymatrix{
    {H_0(\Omega^\infty\Th(-U))} \ar[r]^-{\Omega \omega_*} &
    {H_0(Q(B\Spin(2)_+))} \ar[r]^-{\Omega\trf_*} &
    {H_0(\Omega Q(\RP^\infty_+)) \tensor H_*(\RP^\infty)}
  }
\end{equation*}
which is exact.

Now~\eqref{eq:21} corresponds to the following modified version
of~\eqref{eq:1}
\begin{equation*}
  \xymatrix{
    {\Omega_0^\infty\Th(-U_{\Spin(2)})} \ar[r]^-{\Omega\omega} & {Q_0(B\Spin(2)_+)}
    \ar[r]^-{\Omega_0\partial} & {\Hat\Omega Q(\RP^\infty_+)}.
  }
\end{equation*}
To this fibration there is an associated Eilenberg-Moore spectral sequence
\begin{align}
  \begin{split}
    E^2 &= \Cotor^{H_*(\Hat \Omega_0 Q\RP^\infty_+)}(H_*(Q_0\CP^\infty_+),\F_2)
    \\
    &\isom \Cotor^{H_*(\Hat\Omega_0
      Q\RP^\infty_+)\moddd\Omega\trf_*}(\F_2,\F_2)\tensor H_*(Q(B\Spin(2)_+))
    \modd\Omega\trf_*\\
    &\Rightarrow H_* \Omega^\infty_0\Th(-U)
  \end{split}    \label{eq:9}
\end{align}

\begin{lem}\label{lem:4.1}
  The dual algebra $H^*(\Hat\Omega_0Q\RP^\infty_+) \modd \Omega\trf^*$ is
  polynomial.
\end{lem}
\begin{proof}
  It is a subalgebra of $H^*(\Hat\Omega_0 Q(\RP^\infty_+))$ which
  again is a subalgebra of $H^*(\Omega_0 Q(\RP^\infty_+))$ because
  $\Omega_0 Q(\RP^\infty_+) \simeq \RP^\infty \times \Hat\Omega_0
  Q(\RP^\infty_+)$.  Therefore the lemma follows from
  Theorem~\ref{thm:3}.
\end{proof}
\begin{proof}[Proof of Theorem~\ref{thm:1}]
  From Lemma~\ref{lem:4.1} we get that
  \begin{equation*}
    \Cotor^{H_*(\Hat\Omega_0Q\RP^\infty_+)\moddd \Omega\trf_*}(\F_2, \F_2) =
    E[s^{-1}P(H_*(\Hat\Omega_0Q\RP^\infty)\moddd\Omega\trf_*)]
  \end{equation*}
  Therefore the spectral sequence~\eqref{eq:9} has primitives and
  generators concentrated in $E^2_{0,*}$ and $E^2_{-1,*}$.  Since it
  is a spectral sequence of Hopf algebras, it must collapse.
  Therefore the map
  \begin{equation*}
    \Omega\omega_*\colon  H_*(\Omega_0 \Th(-U)) \to H_*(Q_0(B\Spin(2)_+))
    \modd \Omega\trf_*
  \end{equation*}
  is surjective.
\end{proof}

We next prove Theorem~\ref{thm:5}.  We need a lemma.
\begin{lem}
  \label{lem:4.2}
  The map
  \begin{equation*}
    PH_*(\Tilde\Omega_0 Q\RP^\infty_+) \to P(H_*(\Hat\Omega_0 Q
    \RP^\infty_+)\moddd\Omega\trf_*)
  \end{equation*}
  is surjective.
\end{lem}
\begin{proof}
  First note that $\Hat\Omega_0Q\RP^\infty_+ \simeq \Tilde
  \Omega_0Q\RP^\infty_+\times \RP^\infty$.  The one-dimensional class
  in $PH_*(\RP^\infty)$ is in the image of $P(\Omega\partial_*)$, by
  definition of the double cover $\Tilde \Omega_0Q\RP^\infty_+$, so we
  may substitute $\Hat\Omega_0Q\RP^\infty_+$ for $\Tilde
  \Omega_0Q\RP^\infty_+$ in the statement.  The functor $P$ is left
  exact and has a right derived functor $\Hat P$.  See \cite{G} for a
  survey and references.  The important property is that it vanishes
  when the dual algebra is polynomial.  There is an exact sequence of
  Hopf algebras
  \begin{equation*}
    \xymatrix{
      {\F_2} \ar[r] & {\IM(\Omega\trf_*)} \ar[r] & {H_*(\Hat\Omega_0
        Q\RP^\infty_+)} \ar[r] & {H_*(\Hat\Omega_0 Q\RP^\infty_+)
        \moddd\Omega\trf_*} \ar[r] & {\F_2}
    }
  \end{equation*}
  Now $(\IM(\Omega\trf_*))^* = \IM(\Omega\trf^*)$ is a subalgebra of
  $H^*(Q_0(B\Spin(2)_+))$ and hence is polynomial.  Therefore $\Hat
  P(\IM(\Omega\trf_*)) = 0$ and the lemma follows.
\end{proof}
\begin{cor}
  \label{cor:4.3}
  The map
  \begin{equation*}
    \Cotor^{H_*(\Tilde\Omega_0 Q\RP^\infty_+)}(\F_2, \F_2) \to
    \Cotor^{H_*(\Hat\Omega_0 Q\RP^\infty_+)\moddd\Omega\trf_*}(\F_2, \F_2) 
  \end{equation*}
  is surjective.
\end{cor}
\begin{proof}
  This is because $\Cotor^A(\F_2, \F_2) = E[s^{-1} PA]$ when $A^*$ is
  polynomial.
\end{proof}
\begin{proof}[Proof of Theorem~\ref{thm:5}]
  The spectral sequence gives a filtration $F^0 \supseteq F^{-1}
  \supset \dots$ of $H_*(\Omega_0^\infty\Th(-U))$ which restricts to a
  filtration of $H_*(\Omega_0^\infty\Th(-U))\modd\Omega\omega_*$.  With
  respect to this filtration we have
  \begin{equation*}
    E^0(H_*(\Omega_0^\infty\Th(-U))\modd\Omega\omega_*) \isom
    \Cotor^{H_*(\Hat\Omega_0 Q\RP^\infty_+) \moddd\Omega\trf_*}(\F_2,\F_2)
  \end{equation*}
  There is a map of fibrations
  \begin{equation*}
    \xymatrix{
      {\Omega_0^2 Q\RP^\infty_+} \ar[r]\ar[d] & {*} \ar[r]\ar[d] &
      {\Tilde\Omega_0 Q\RP^\infty_+}\ar[d]\\
      {\Omega_0^\infty\Th(-U)}\ar[r]^-{\Omega\omega} & {Q_0(B\Spin(2)_+)}
      \ar[r]^-{\Omega\trf} & {\Hat\Omega_0 Q\RP^\infty_+}
    }
  \end{equation*}
  and an associated map of spectral sequences which on the $E^2$-term is
  \begin{equation*}
    \Cotor^{H_*(\Tilde \Omega_0 Q\RP^\infty_+)}(\F_2,\F_2) \to 
    \Cotor^{H_*(\Hat\Omega_0 Q\RP^\infty_+)\moddd\Omega\trf_*}(\F_2,\F_2) \tensor
    H_*(Q(B\Spin(2)_+))\modd\trf_*
  \end{equation*}
  Since both spectral sequences collapse, we get that the map
  \begin{equation}
    \label{eq:14}
    H_*(\Omega_0^2 Q\RP^\infty_+) \to H_*(\Omega^\infty_0\Th(-U))\modd\Omega
    \omega_*
  \end{equation}
  is filtered and on filtration quotients the map is identified with 
  \begin{equation*}
    \Cotor^{H_*(\Tilde\Omega_0 Q\RP^\infty_+)}(\F_2,\F_2) \to 
    \Cotor^{H_*(\Hat\Omega_0 Q\RP^\infty_+)\moddd\Omega\trf_*}(\F_2,\F_2)
  \end{equation*}
  Since this is surjective by Lemma~\ref{lem:4.2}, then also the
  map~\eqref{eq:14} is surjective.

  This proves (i).  (ii) is just the dual statement of (i).  To prove
  (iii) we see that the quotient
  \begin{equation*}
    QH_*(\Omega^2 Q(\RP^\infty_+)) \to Q(H_*(\Omega^\infty\Th(-U))
    \modd \Omega\omega_*)
  \end{equation*}
  is identified under suspension with
  \begin{equation*}
    PH_*(\Omega Q( \RP^\infty_+)) \to P(H_*(\Omega
    Q(\RP^\infty_+))\moddd \Omega\trf_*)
  \end{equation*}
  which again by suspension is mapped to
  \begin{equation*}
    PH_*(Q(\RP^\infty_+)) \to P(H_*(Q(\RP^\infty_+))\moddd\trf_*) =
    \Cok(P\trf_*).
  \end{equation*}
  By dualising we get $\Ker(Q\trf_*)$ as claimed.
\end{proof}

\section{Proof of Theorem \ref{thm:2}}
\label{sec:proof-theorem}

We know from theorem~\ref{thm:3} that $H_*(\Omega Q(\RP^\infty_+))$
is an exterior algebra.  We also know that $H_*(QB\Spin(2)_+)$ is a
polynomial algebra.  We have the following commutative diagram
\begin{equation*}
  \xymatrix{
    {QH_*(Q(B\Spin(2)_+))} \ar[r]^-{Q(\Omega \trf_*)} \ar[d]^{\cong}
    & {QH_*(\Omega Q(\RP^\infty_+))} \ar[d]^{\cong}\\
    {PH_*(Q\Sigma(B\Spin(2)_+))} \ar[r]^-{P(\trf_*)} 
    & {PH_*(Q(\RP^\infty_+))} \\
  }
\end{equation*}
It follows that $Q(\Omega\trf_*)$ is injective.  These three facts
prove that $H_*(Q\CP^\infty_+)\modd \Omega\trf_* = \xi
H_*(Q\CP^\infty_+)$.  We calculate the Becker-Gottlieb transfer of the
bundle
\begin{equation*}
  E\Spin(2) \times_{\Spin(2)} S^2 \to B\Spin(2).
\end{equation*}
\begin{lem}
  Let $N,S\colon B\Spin(2) \to E\Spin(2)\times_{Spin(2)} S^2$ denote the
  sections at the north and south pole, respectively.  Then the
  Becker-Gottlieb transfer is
  \begin{equation*}
    \tau = N+S \in [B\Spin(2), Q(E\Spin(2)\times_{\Spin(2)} S^2 _+)]
  \end{equation*}
\end{lem}
\begin{proof}
  This is similar to the Becker-Gottlieb calculations in \cite{GMT}:
  $S^2$ is the $\Spin(2)$-equivariant pushout of $D^2 \leftarrow S^1
  \rightarrow D^2$ and therefore the bundle $E\Spin(2)
  \times_{\Spin(2)} S^2$ is the fibrewise pushout of $E\Spin(2)
  \times_{\Spin(2)} D^2 \leftarrow E\Spin(2) \times_{\Spin(2)} S^1
  \rightarrow E\Spin(2) \times_{\Spin(2)} D^2$.  Then properties
  (A1)--(A3) in \cite[p.\ 15]{GMT} proves the proposition.  Indeed the
  transfer of $E\Spin(2) \times_{\Spin(2)} S^1$ vanishes by (A3) and
  the transfer of $E\Spin(2) \times_{\Spin(2)} D^2$ is the section at
  the center of $D^2$ by (A1).  Then the additivity (A2) proves that
  the transfer of the whole bundle is $N+S$.
\end{proof}
\begin{cor}
  Let $E\Spin(2)\times_{\Spin(2)} S^2 \to B\Spin(2)$ classify the
  vertical tangentbundle.  Then 
  \begin{equation*}
    \alpha = \iota + c \in[B\Spin(2), Q(B\Spin(2)_+)]
  \end{equation*}
  where $\iota$ is the usual inclusion of $B\Spin(2)$ and $c$ is the
  orientation reversal map.
\end{cor}
\begin{proof}[Proof of theorem~\ref{thm:2}]
  We have $\iota_* a_i = a_i$ and $c_* a_i = (-1)^i a_i$.  Therefore
  \begin{equation*}
    (\iota + c)_* a_i = \sum_{r+s=i} (-1)^s a_ra_s
  \end{equation*}
  Reducing mod 2 we get
  \begin{equation*}
    (i+c)_*(a_{2i}) = a_i^2 \quad\text{and}\quad (i+c)_*a_{2i+1} = 0
  \end{equation*}
  Since $B\Spin(2)\to B\Spin(3)$ maps $a_{2i} \mapsto b_i$, we have
  proved that the composition in theorem~\ref{thm:2} maps $b_i$ to
  $a_i^2$ as claimed.

  Then it will also map $Q^{2I} b_i$ to $(Q^I a_i)^2$ and hence the
  composition is surjective.

  Both $H_*(Q(B\Spin(2)_+))$ and $H^*(Q_0(B\Spin(2)_+))$ are free
  commutative.  This follows from the fact that $\lambda\colon 
  H_*(B\Spin(2)) \to H_*(B\Spin(2))$ is surjective, similarly to the
  case of $Q(\RP^\infty_+)$.  But then
  \begin{equation*}
    \xi\colon  H_*(Q(B\Spin(2)_+)) \to \xi H_*(Q(B\Spin(2)_+))
  \end{equation*}
  is an isomorphism so the same holds for $\xi H_*(Q(B\Spin(2)_+))$.
\end{proof}
\begin{proof}[Proof of Corollary~\ref{cor:1.8}]
  By the exact sequence in Theorem~\ref{thm:1} and by
  Theorem~\ref{thm:5}, the kernel of
  \begin{equation*}
    H^*(\Omega^\infty_0\Th(-U)) \to H^*(\Omega^2_0 Q(\RP^\infty_+))
  \end{equation*}
  is exactly $H^*(Q_0(B\Spin(2)_+)) \moddd \Omega\trf^*$.  By
  theorem~\ref{thm:2}~(ii), this injects into $H^*(Q(B\Spin(3)_+))$.
\end{proof}

\newpage

\section{Adapting \cite{MW}}
\label{sec:adapting-mw}

This is the second part of the paper, and the aim is to adapt the
proof in \cite{MW}.  As explained in the introduction we can let $\theta =
\theta_{\SO}$ and then for genus $\geq 2$ we have $\MM^\theta(F,\gamma)
= B\Gamma(F,\gamma)$, where $\gamma$ is an orientation of $F$ and
$\Gamma(F,\gamma) = \pi_0 \Diff(F,\gamma)$ is the oriented mapping
class group of $F$.  Then we can let $F = F_{g,2}$ and let
$\Gamma_{\infty,2} = \colim \Gamma_{g,2}$ where the colimit is over
glueing an oriented torus.  Then \cite{MW} proves that there is a homology
equivalence
\begin{equation*}
  \Z\times B\Gamma_{\infty,2} \to \Omega^\infty\Th(-U_{SO})
\end{equation*}

For $\theta = \theta_\Spin$ we can again let $F = F_{g,2}$ and let
\begin{equation*}
  \MM^\theta(F_{\infty,2}) := \hocolim \MM^\theta(F_{g,2})
\end{equation*}
where the hocolim is over glueing a torus.  There are two essentially
different ways of doing this because we can choose either an Arf
invariant 0 torus or an Arf invariant 1 torus.  Which one we use is
not important however, because the composition of two tori will be a
surface of genus 2 and with an Arf invariant 0 spin structure anyhow.

Then we adapt the proof to showing that there is a homology
equivalence
\begin{equation*}
  \Z\times \MM^\theta(F_{\infty,2}) \to
  \Omega^\infty\Th(-U_{\Spin(2)})
\end{equation*}
Since the $\MM^\theta(-)$ satisfies Harer stability, we also get that
$\MM^\theta(F,\gamma)$ has the same homology as
$\Omega^\infty\Th(-U_{\Spin(2)})$ in a stable range.

Most of the modifications are straightforward and the proofs are valid
for any vectorbundle $\theta\colon  U_3 \to B_3$.  Only at the very end
shall we specialise to the case $\theta = \theta_\Spin$.  The idea is
roughly as follows.  All the sheaves in \cite{MW} are made out of either
submersions $\pi\colon  E \to X$ with oriented three-dimensional fibres, or
surface bundles $q\colon  M \to X$ with oriented two-dimensional fibres,
with some extra structure.  Then we can modify the definition by
removing the word ``oriented'' anbd instead include a bundle map
$T^\pi E \to U_3$ or $T^q M \to U_2$.  The original case in \cite{MW} can
the be recovered by setting $\theta = \theta_{\SO}$ (the sheaves will
be slightly fattened versions of those in \cite{MW}).

This procedure works very well, and for most of the chapters we shal
just give the modified definitions and claim that the proofs work in
our more general situation as well.  There is one point that needs
attention, however.  Namely the definition of $\Erg$ and the sheaf map
$\LT \to \WT$ in \cite[Chapter 5]{MW}.  To do this properly in the
added generality we shall need to give a new definition of fibrewise
surgery.  Also \cite[Chapter 6]{MW} about the ``connectivity problem''
need some attention.

\section{The sheaves}
This section defines the appropriate generalisations of the sheaves on
$\mathscr{X}$ defined in \cite[Section 2]{MW}.  Let $\theta\colon U_3
\to B_3$ be a 3-dimensional real vectorbundle.  Let $q_0 \colon T(S^1
\times [0,1] \times \R) \to U_3$ be a fixed bundle map, constant in
the $[0,1]\times \R$-directions.
\begin{defn}
  \label{defn:sheaves}
  Let $\mathcal{V}^\theta$ be the sheaf on $\mathscr{X}$ defined such
  that $\mathcal{V}(X)$ is the set of $(\pi,f,q)$ such that $(\pi,f)\colon 
  E^{k+3} \to X^k \times\R$ is a proper smooth map, $\pi\colon  E \to X$ is
  a graphic submersion, $f$ is fibrewise regular, and $q\colon  T_\pi E \to
  U_3$ is a bundle map.  We assume that near the boundary of $E$,
  $(\pi,f)$ agrees over $X \times \R$ with $S^1 \times [0,1] \times
  \R$ and $q$ agrees with $q_0$.
  
  Define $h\mathcal{V}^\theta$, $\mathcal{W}^\theta$,
  $h\mathcal{W}^\theta$, $\mathcal{W}_\mathrm{loc}^\theta$ and
  $h\mathcal{W}_\mathrm{loc}^\theta$ similarly.
\end{defn}

\begin{rem}
  For $\theta = E\SO(3) \times_{\SO(3)} \R^3 \to B\SO(3)$, the map $q\colon 
  T_\pi E \to U$ induces an orientation on the fibres of $\pi\colon  E \to
  X$.  Thus for this $\theta$ there is a sheaf map
  \begin{equation*}
    \mathcal{V}^\theta \to \mathcal{V}
  \end{equation*}
  and this is a weak equivalence.  Thus $\mathcal{V}^\theta$ is a fat
  version of the sheaf $\mathcal{V}$ in \cite{MW}.
\end{rem}

Following \cite{MW} we get a diagram of classifying spaces
\begin{equation*}
  \xymatrix{
    {|\mathcal{V}_c^\theta|} \ar[r]\ar[d]& {|\mathcal{W}^\theta|}
    \ar[r]\ar[d] &
    {|\mathcal{W}_{\mathrm{loc}}^\theta|}\ar[d]\\
    {|h\mathcal{V}^\theta|} \ar[r]& {|h\mathcal{W}^\theta|} \ar[r]&
    {|h\mathcal{W}_{\mathrm{loc}}^\theta|}
  }
\end{equation*}
where the vertical maps are induced by taking the 2-jet prolongation
of $f$.

We aim at generalising \cite{MW} to the statement $\Omega B
|\mathcal{V}^\theta_c| \simeq |h\mathcal{V}^\theta|$.

\begin{lem}
  Let $\Th(-U_2)$ denote the Thom spectrum of the vitual bundle $-U_2$
  over $B(2)$.  Then
  \begin{equation*}
    |h\mathcal{V}^\theta| \simeq \Omega^\infty\Th(-U_2)
  \end{equation*}
\end{lem}
\begin{lem} We have
  \begin{equation*}
    |\mathcal{V}_c^\theta| \simeq \coprod_{F}
     \MM^\theta(F).
  \end{equation*}
  where the disjoint union is taken over surfaces with two boundary
  components, one in each diffeomorphism class.
  
  If $U\to B$ is orientable, this means that the disjoint union is
  over the surfaces $F_{g,2}$, $g \geq 0$.
\end{lem}
\begin{proof}
  This is proved similarly to the case considered in \cite{MW}.
\end{proof}

\section{Adjusting the proof}
Most of the proof given in \cite{MW} goes through with little or no change
also in this more general situation.  We describe the necessary
changes chapter for chapter.

\subsection{Chapter 3}
\cite{MW} determines the homotopy types of $|h\mathcal{V}|$,
$|h\mathcal{W}|$ and $|h\mathcal{W}_\mathrm{loc}|$ and proves that 
\begin{equation*}
  |h\mathcal{V}| \to
  |h\mathcal{W}| \to |h\mathcal{W}_\mathrm{loc}|
\end{equation*}
is a homotopy fibre sequence.

Let $\Bun(\R^3,U_3)$ denote the space of bundle maps from $\R^3$,
considered as a bundle over a point, to the bundle $U_3$.  As in
\cite{MW} we let $S(\R^3)$ be the vectorspace of quadratic forms on
$\R^3$ and $\Delta\subseteq S(\R^3)$ be the subset of degenerate
quadratic forms.

Define an $\OO(3)$-space $A^\theta(\R^3)$ by
\begin{equation*}
  A^\theta = ((\R^3)^* \times S(\R^3) - \{0\} \times \Delta) \times
  \Bun(\R^3, U_3)
\end{equation*}
and define
\begin{equation*}
  \GW^\theta (3,n)) = \OO(n+3) \times_{\OO(n)\times\OO(3)}
  A^\theta (\R^3).
\end{equation*}
Thus a point in $\GW(3,n)$ is a quadruple $(V,l,q,\xi)$ where
$V\subseteq \R^{3+n}$ is a three-dimensional subspace, $l\colon V \to
\R$ is a linear map, $q\colon V \to \R$ is a quadratic map, and $\xi\colon 
V \to U_3$ is a bundle map, subject to the condition that $q$ is
non-degenerate if $l = 0$.

\begin{example}
  For $U_3 = E\OO(3)\times_{\OO(3)} \R^3$, the space $A^\theta(\R^3)$
  has the same (equivariant) homotopy type as the $A(\R^3)$ of
  \cite{MW}.  For $U_3 = E\SO(3) \times_{\SO(3)} \R^3$, the space
  $\GW^\theta (3,n)$ has the same homotopy type as $\GW(3,n)$ of
  \cite{MW}.
\end{example}

Let $\Sigma^\theta(3,n) \subseteq  \GW^\theta (3,n)$ be the
subspace corresponding to $\{0\} \times (S(\R^3) - \Delta) \times
\Bun(\R^3, U_3)\subseteq A^\theta(\R^3)$, and let
\begin{equation*}
  \GV^\theta(3,n) = \GW^\theta(3,n) - \Sigma^\theta(3,n)
\end{equation*}

Let $\U^\theta_{n} \to \GW^\theta(3,n) $ be the universal
bundle.  We get a cofibration sequence
\begin{equation*}
  \thh({\U_{n}^\theta}^\perp \vert \GV^\theta (3,n)) \to 
  \thh({\U_{n}^\theta}^\perp) \to 
  \thh({\U_{n}^\theta}^\perp \oplus {\U_{n}^\theta}^* \vert \Sigma^\theta(3,n) ) 
\end{equation*}
and an associated fibration sequence of infinite loop spaces
\begin{equation*}
  \Omega^\infty \spectrumfont{h\V}^\theta \to
  \Omega^\infty \spectrumfont{h\W}^\theta \to
  \Omega^\infty \spectrumfont{h\Wloc}^\theta
\end{equation*}
as in \cite[Paragraph 3.1]{MW}.

We have the following generalisations of \cite{MW}:
\begin{thm}\ 
  \begin{enumerate}[(i)]
  \item $|h\W^\theta| \simeq \Omega^\infty\spectrumfont{h\W}^\theta$
  \item $|h\V^\theta| \simeq \Omega^\infty\spectrumfont{h\V}^\theta$
  \item $|h\Wloc^\theta| \simeq \Omega^\infty\spectrumfont{h\Wloc}^\theta$
  \item $|\Wloc^\theta| \simeq \Omega^\infty\spectrumfont{h\Wloc}^\theta$
  \end{enumerate}
\end{thm}
\begin{proof}
  Similar to \cite{MW}.
\end{proof}

\subsection{Chapter 4}

In 4.2, we define $\W_\theta^\AA$ and $h\W_\theta^\AA$ in the obvious
way.  These are sheaves of posets.

In 4.3, we define a sheaf $\mathcal{T}^\AA_\theta$ as in
\cite[Definition 4.3.1]{MW}, but with the added data of a bundle map
$q\colon  T^\pi E \to U_3$.  Notice that this is a small errata to
\cite{MW}: Their $\mathcal{T}^\AA$ should consist of $(\pi, \psi)\colon  E
\to X\times\R$ such that $\pi\colon  E \to X$ is a submersion with
\emph{oriented} fibres.

With these modifications, the proof in \cite[Section 4.3]{MW} goes
through without further difficulties.  Thus we get
\begin{thm}
  $|\W^\theta| \simeq |h\W^\theta|$
\end{thm}

\subsection{Chapter 5:  Surgery}
\cite[Chapter 5.2]{MW} is about fibrewise surgery.  The idea is
roughly as follows.  Given a bundle $q\colon  M \to X$ of manifolds, a
finite set $T$, a Riemannian vectorbundle $\omega\colon  V \to T \times X$
with isometric involution $\rho\colon  V \to V$, and an embedding $e\colon 
D(V^\rho) \times_{T\times X} S(V^{-\rho}) \to M - \partial M$, then
one performs surgery by removing the interior of the embedded
$D(V^\rho) \times_{T\times X} S(V^{-\rho})$ and replacing it with
$S(V^\rho) \times_{T\times X} D(V^{-\rho})$.

In our generalised setting, $M$ will be equipped with a bundle map
$\xi\colon  T^q M \to U_2$.  We would like to perform surgery in a way that
we end up with a bundle $\bar q\colon  \bar M \to X$, equipped with a bundle
map $\bar \xi\colon  T^{\bar q} \bar M \to U_2$.  We describe how to do
this.

\subsubsection{Saddles}
\label{sec:saddles}

Choose once and for all a smooth function $\tau\colon  [0,1] \to [0,1]$
which is 0 near 0 and 1 near 1.  Let $Y$ be a manifold and $\omega\colon  V
\to Y$ a Riemannian vectorbundle with isometric involution $\rho\colon  V
\to V$.  Let $g\colon  Y \to \R$ be smooth.  As in \cite{MW} we define the saddle
of $V$ to be the subset
\begin{equation*}
  \Sad(V) = \{ v \in V \mid |v_+||v_-| \leq 1\}  
\end{equation*}
Define three smooth functions by
\begin{align*}
  f_0(v) & = g\omega(v) + |v_+|^2 - |v_-|^2\\
  f_+(v) & = g\omega(v) + \frac1{|v_-|^2} \left( |v_+|^2 |v_-|^2
    \tau(|v_+||v_-|)
    + (1-\tau(|v_+||v_-|))\right) - |v_-|^2\\
  f_-(v) & = g\omega(v) + |v_+|^2 - \frac1{|v_+|^2} \left( |v_+|^2
    |v_-|^2 \tau(|v_+||v_-|)
    + (1-\tau(|v_+||v_-|))\right)\\
\end{align*}

The map $f_0$ is defined on all of $\Sad(V)$ and is fibrewise regular
except at the zero section of $V$, where it has a Morse singularity
with critical value given by $g\omega$.  The maps $f_{\pm}$ is defined
on $\Sad(V) - V^{\pm \rho}$, is fibrewise regular and proper, and
agrees with $f_0$ near $\partial \Sad(V)$.  The following picture
shows the level curves of $f_0$ in $V$.  $\Sad(V)\subseteq V$ is the
shaded area.
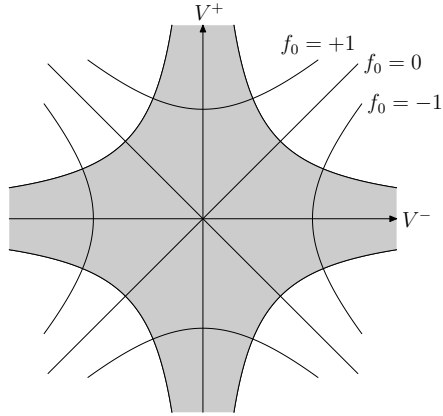
\begin{center}
  \includegraphics[width=0.4\linewidth]{levelcurves.1}  
\end{center}
This should be compared with the level curves of $f_+$ and $f_-$,
shown in the following pictures.
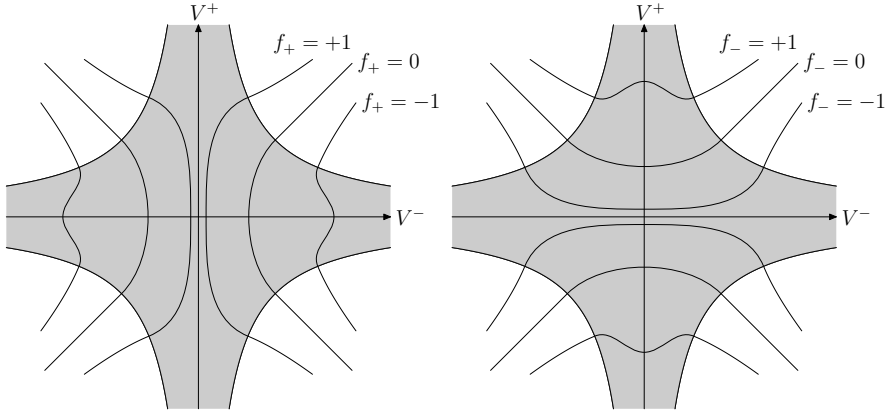
\begin{center}
  \includegraphics[width=0.4\linewidth]{levelcurves.2}  
  \includegraphics[width=0.4\linewidth]{levelcurves.3}  
\end{center}

  Moreover, $f_+$ defines a
diffeomorphism
\begin{gather}
  \label{eq:15}
  \begin{split}
    \Sad(V) - V^\rho & \to D(V^\rho)\times_Y S(V^{-\rho}) \times \R\\
    v & \mapsto (|v_-|v_+, |v_-|^{-1}v_-, f_+(v))
  \end{split}
\end{gather}
Similarly, $f_-$ defines a diffeomorphism
\begin{gather}
  \label{eq:16}
  \begin{split}
    \Sad(V) - V^{-\rho} & \to S(V^\rho)\times_Y D(V^{-\rho}) \times \R\\
    v & \mapsto (|v_+|v_+, |v_+|^{-1}v_-, f_-(v))
  \end{split}
\end{gather}

\begin{rem}
  Comparing \eqref{eq:15} to equation \cite[equation (5.3)]{MW} we see that,
  up to diffeomorphism, the process of removing $V^\rho$ and replacing
  $f$ with $f_+$ is equivalent to glueing $D(V^\rho) \times _Y
  S(V^{-\rho}) \times \R$ to $\Sad(V) - V^\rho$ along \cite[equation
  (5.3)]{MW}.  Similarly for~\eqref{eq:16} and \cite[equation
  (5.4)]{MW}.
\end{rem}

\begin{defn}
  Given a vectorbundle $\omega\colon  V \to Y$ and a smooth $g\colon  Y \to \R$ as
  above, we let
  \begin{align*}
    M_+(V,g) &= f_+^{-1}(0) \subseteq \Sad(V)\\
    M_-(V,g) &= f_-^{-1}(0) \subseteq \Sad(V)\\
  \end{align*}
\end{defn}
By our earlier remarks we see that both $M_+(V,g)$ and $M_-(V,g)$ agrees
near $\partial \Sad(V)$ with $f_0^{-1}(0)$.  By restriction
of~\eqref{eq:15} we get a diffeomorphism over $Y$
\begin{align*}
  M_+(V,g) \to D(V^\rho) \times_Y S(V^{-\rho})
\end{align*}
and the fibrewise differential induces an isomorphism
\begin{align*}
  T^\omega M_+(V,g) \times \R \to T^\omega V_{|M_+(V,g)}
\end{align*}
Similarly for $M_-(V,g)$.

This gives an alternative description of surgery.  Namely, given a
surface bundle $q\colon  M \to X$, a finite set $T$, a Riemannian
vectorbundle $V \to T \times X$ with isometric involution $\rho\colon  V \to
V$, a smooth function $g\colon  T\times X \to \R$, and an embedding over $X$
$\lambda\colon  M_+(V,g) \to M - \partial M$, then one performs surgery by
replacing the embedded $M_+(V,g)$ by $M_-(V,g)$.  Since $M_+(V,g)$ and
$M_-(V,g)$ agree near their boundary, this gives a welldefined smooth
bundle $\bar q\colon  \bar M \to X$.  Moreover the following is true.  If
$M$ is equipped with a bundle map $ \xi\colon  T^q M \to U_2$ and $V$ is
equipped with $\xi\colon  T^\omega V_{|\Sad(V)} \to U_3$, and the fibrewise
differential of $\lambda$ is over $U_2$, then $\bar M$ gets a
canonical map $T^{\bar q}\bar M \to U_2$.

\subsubsection{The sheaves}
\label{sec:sheaves}

Keeping these remarks in mind, we make the following definitions.
$\WlocT^\theta(X)$ is the set of
\begin{enumerate}[(i)]
\item $\omega\colon  V \to T \times X$ a Riemannian vector bundle with
  isometric involution $\rho$, as in \cite{MW}.
\item $g\colon  T \times X \to \R$ a smooth function.
\item $\xi\colon  T^\omega V_{|\Sad(V)} \to U_3$ a vectorbundle map
\end{enumerate}
and $\WT^\theta(X)$ is the set of
\begin{enumerate}
\item $(V,g,\xi) \in \WlocT(X)$
\item $q\colon  M \to X$ a bundle of surfaces
\item $\xi\colon  T^q M \to U_2$ a bundle map
\item $e\colon  M_+(V,g) \to M - \partial M$ an embedding over $X$ such that
  the fibrewise differential $De$ is over $U_2$.
\end{enumerate}

\subsubsection{The proofs}
\label{sec:proofs}
We go through the definitions and proofs in \cite[Chapter 5]{MW} and
describe what modifications are needed in this more general
situation.  Again this is summarised in the diagram
\begin{gather}
  \label{eq:17}
  \begin{split}
  \xymatrix{
    {\W^\theta} \ar[rr] & & {\Wloc^\theta} \\
    {\LL^\theta} \ar[rr]\ar[u] & & {\Lloc^\theta}\ar[u]\\
    {\hocolim \LT^\theta} \ar[u]\ar[d]\ar[rr]& & {\hocolim \LlocT^\theta}\ar[u]\ar[d]\\
    {\hocolim \WT^\theta} \ar[rr]&& {\hocolim \WlocT^\theta}
  }
  \end{split}
\end{gather}

\subsubsection{Second row}
\label{sec:second-row}

Define $\Lloc^\theta(X)$ to be the set of
\begin{enumerate}[(i)]
\item $(p,g)\colon  Y \to X\times \R$ proper smooth maps such that $p$ is
  etale and graphic and such that $g$ is smooth.
\item $\omega\colon  V \to Y$ is a Riemannian vectorbundle with isometric
  involution $\rho\colon  V \to V$.
\item $\xi\colon  T^\omega V _{|\Sad(V)} \to U_3$ a bundle map
\end{enumerate}
and let $\LL^\theta(X)$ be the set of
\begin{enumerate}[(i)]
\item $(\pi,f,\xi) \in \W^\theta(X)$ with $(\pi,f)\colon  E \to X \times \R$,
  $\xi\colon  T^\pi E \to U_3$.
\item $(p,g,V,\xi) \in \Lloc^\theta(X)$
\item $\lambda\colon  \Sad(V) \to E - \partial E$ an embedding over $X
  \times \R$ such that the fibrewise differential $D\lambda$ is over
  $U_3$.
\end{enumerate}
The proofs given in \cite{MW} of \cite[Proposition 5.3.3]{MW} and
\cite[Proposition 5.3.7]{MW} goes through with the obvious changes and
proves that the sheaf maps $\Lloc^\theta \to \Wloc^\theta$ and
$\LL^\theta \to \W^\theta$ are weak equivalences.

\subsubsection{Third row}
\label{sec:third-row}

Let $\LlocT^\theta(X)$ be the set of
\begin{enumerate}[(i)]
\item $(p,g,V,\xi) \in \Lloc^\theta(X)$
\item $h\colon  S \times X \to Y$ an embedding over $3 \times X$
\item $\delta\colon  Y - \IM(h) \to \{\pm 1\}$ continuous
\end{enumerate}
and let $\LT^\theta(X)$ be the set of
\begin{enumerate}[(i)]
\item $(p,g,V,\xi,h,\delta) \in \LlocT^\theta(X)$
\item $(\pi,f,\xi) \in \W^\theta(X)$
\item $\lambda\colon  \Sad(V) \to E - \partial E$ embedding over $X \times
  \R$ such that the fibrewise differential $D\lambda$ is over $U_3$.
\end{enumerate}
The proofs given in \cite{MW} of \cite[Proposition 5.4.2]{MW} and
\cite[Proposition 5.4.4]{MW} goes through with the obvious changes and
proves that the sheaf maps $\hocolim \LlocT^\theta \to \Lloc^\theta$
and $\hocolim \LT^\theta \to \LL^\theta$ are weak equivalences.

\subsubsection{Fourth row, right hand column}
\label{sec:fourth-row-right}

\cite[Lemma 5.5.2]{MW} and \cite[Corollary 5.5.3]{MW} goes through as
in \cite{MW}.

\subsubsection{Fourth row left hand column}
\label{sec:fourth-row-left}

This is more technical, and more changes are needed to adapt the proof
in \cite{MW}.  The modified definitions of $\WT^\theta$ and
$\WlocT^\theta$ were made with this in mind.  The problem is to give a
definition of $\Erg$ and to define a sheaf map $\LT^\theta \to
\WT^\theta$, natural in $T \in \KK$.

Take an element of $\LT^\theta(X)$.  This consists of
$(p,g,V,\xi,h,\delta) \in \LlocT^\theta(X)$, $(\pi,f,\xi) \in
\W^\theta(X)$, and $\lambda\colon  \Sad(V) \to E - \partial E$.  Define
$Y_0, Y_+, Y_- \subseteq Y$ and $V_+, V_-, V_0\subseteq V$ as in \cite{MW}.
Define $\Erg$, $\frg$ in the following way
\begin{itemize}
\item On the embedded $\Sad(V_+)$, remove $V_+^\rho$ and replace $f$
  by $f_+$.
\item On the embedded $\Sad(V_-)$, remove $V_-^{-\rho}$ and replace
  $f$ by $f_-$.
\item On the embedded $\Sad(V_0)$, remove $V_+^\rho$ and replace $f$
  by $f_+$.
\end{itemize}
This defines a bundle $(\pi^\mathrm{rg},\frg)\colon  \Erg \to X \times\R$ of
smooth compact surfaces.  Now let $M = (\frg)^{-1}(0)$.  This is a
bundle of smooth compact surfaces over $X$, and is equipped with the
following extra structure
\begin{enumerate}[(i)]
\item A bundle map $\xi\colon  T^\pi M \to U_2$
\item A Riemannian vectorbundle $\omega\colon  h^* V_0 \to T \times X$ with
  isometric involution $\rho$.
\item A bundle map $\xi\colon  T^\omega(h^*V_0)_{|\Sad(V)} \to U_3$.
\item A smooth function $g\colon  T \times X \to Y \to \R$.
\item An embedding (over $X$) $e\colon  M_+(V,g) \to M$ such that the fibrewise
  differential is over $U_2$.
\end{enumerate}
That is, we have an element of $\WT^\theta(X)$.  This defines a sheaf
map $\LT^\theta \to \WT^\theta$ which is natural in $T \in \KK$.  Just
as in \cite{MW} one proves that $\LT^\theta \to \WT^\theta$ is an
equivalence.

\subsubsection{Using the concordance lifting property}
\label{sec:using-conc-lift}

To prove that the sheaf maps $\WT^\theta \to \WlocT^\theta$ has the
concordance lifting property we need the following lemma
\begin{lem}
  Let $A\subseteq X$ be a cofibration and let $V \to [0,1]\times X$ be
  a vectorbundle.  Let $U \to B$ be another vectorbundle.  Then any
  bundle map $\xi\colon  V_{|\{0\}\times X\cup [0,1]\times X} \to U$ extends
  to a bundle map $V \to U$
\end{lem}
\begin{proof}
  Choose a retraction $r\colon  [0,1]\times X \to \{0\} \times X \cup
  [0,1]\times A$.  Now the fibre bundle $\mathrm{Iso}(V,r^* V) \to
  [0,1]\times A$ has a canonical section over $\{0\}\times X \cup
  [0,1]\times A$.  This section extends over all of $[0,1]\times X$
  because $\{0\}\times X\cup [0,1]\times A \to [0,1]\times X$ is a
  trivial cofibration.  This section defines a bundle map
  \begin{equation*}
    \xymatrix{
      {V} \ar[r]^-{\hat r} \ar[d] & {V_{|\{0\}\times X\cup [0,1]\times
          A}}\ar[d]\\
      {[0,1]\times X} \ar[r]^-{r} & {\{0\}\times X\cup [0,1]\times A}
    }
  \end{equation*}
  and we can compose $\xi$ with $\hat r$.
\end{proof}
\begin{prop}
  The map $\WT^\theta \to \WlocT^\theta$ has the concordance lifting
  property.
\end{prop}
\begin{proof}
  Let $\chi\in\WT^\theta$ be an element given by
  \begin{itemize}
  \item $(V,\xi)\in\WlocT^\theta(X)$ with $\omega\colon  V \to T\times X$
    and $\xi\colon  T^\omega V_{|\Sad(V)} \to U_3$
  \item $q\colon  M \to X$ a surface bundle (with certain boundary
    conditions).
  \item $\xi\colon  T^q M \to U_2$ a bundle map
  \item $e\colon M_+(V,g) \to M$ an embedding over $X$ such that the
    fibrewise differential $De$ is over $U_2$.
  \end{itemize}
  Suppose given a concordance of $(V,\xi)$.  This will be given by a
  vectorbundle $\Tilde \omega\colon  \Tilde V \to (0,1)\times T\times X$ and
  $\Tilde \xi\colon  T^\omega \Tilde V_{|\Sad(\Tilde V)} \to U$.  We can
  choose an isomorphism $\Tilde V \cong (0,1)\times V$ over
  $(0,1)\times T \times X$.  Put $\Tilde M = (0,1)\times M$ and
  $\Tilde q = (0,1)\times q$.  Let $\Tilde g = g\circ
  \mathrm{pr}_{T\times X}\colon  (0,1)\times T\times X \to \R$.  Then we
  have the isomorphism $M_+(\Tilde V, \Tilde g) \cong
  (0,1)\times M_+(V,g)$ and we can set $\Tilde e = (0,1)\times e\colon 
  M_+(\Tilde V, \Tilde g) \cong (0,1)\times M_+(V,g) \to (0,1)\times
  M$.

  It remains to define a bundle map $T^{\Tilde q}\Tilde M \to U_2$
  which is specified on $T^{\Tilde q}\Tilde M_{|\{0\}\times M\cup
  [0,1]\times M_+(V,g)}$.  This can be done by the previous lemma,
  using that $M_+(V,g) \to M$ is a cofibration.
\end{proof}

\subsection{Chapter 6: The connectivity problem}

We describe how to adapt the definition of the sheaf $\CM$ and prove
that it is contractible.  Let $\R^2\times\R$ have the standard
euclidean metric and involution $\rho = \mathrm{diag}(1,1,-1)$.  For
any finite set $T$ and a manifold $X$ we have the trivial vectorbundle
$V = \R^2\times\R\times T \times X$ over $T\times X$ and we have
canonical identifications
\begin{itemize}
\item $\Sad(V) = \Sad(\R^2\times\R) \times T \times X$.
\item $T^\omega V _{|\Sad(V)} = \R^2 \times\R\times\Sad(V)$
\item $D^2 \times S^0 \times T \times X \cong M_+(0) \subseteq \Sad(V)$
\end{itemize}
Thus to promote $V$ to an element of $\WlocT^\theta$ with $T\to \{1\}$
we must specify a bundle map $T^\omega V _{|\Sad(V)} \to U_3$, or
equivalently a map $\Sad(V) \to \Bun(\R^2\times\R,U_3)$.

\begin{defn}
  Let $M$ be a surface and $TM \to U_2$ a bundle map.  Let $\CMt$ be the
  sheaf whose value at a connected manifold $X$ is the set of
  \begin{itemize}
  \item A finite set $T$
  \item A map $\Sad(V) \to \Bun(\R^2\times\R,U_3)$, where $V = \R^2
    \times\R\times T\times X$ as above
  \item An embedding $e_T\colon  M_+(0) \to (M-\partial M) \times X$ over
  $X$ such that the fibrewise differential $De_T$ is over $U_2$ and
  such that surgery along $e_T$ results in a connected surface bundle
  over $X$.
  \end{itemize}
\end{defn}

We want to prove that $B|\CMt| \simeq |\beta \CMt^\mathrm{op}|$ is
contractible.  We proceed as in \cite{MW}: Given a closed set $A \subseteq
X$ and a germ $s_0 \in \colim_U \beta\CMtop(U)$ we extend this germ to
an element of $\beta\CMtop(X)$.  The germ $s_0$ consists of a locally
finite open cover $(U_j)_{j\in J}$ of $U$ and objects $\phi_{RR} \in
\CMt(U_R)$ for each finite non-empty $R\subseteq J$, and for each
$R\subseteq S$ a morphism $\phi_{RS}\colon  \phi_{SS} \to
{\phi_{RR}}_{|U_S}$ satisfying the cocycle condition.  Each of the
$\phi_{RR}$ defines an embedding
\begin{equation*}
  D^2 \times S^0 \times T_R \times U_R \cong M_+(0) \to (M-\partial M)
  \times U_R
\end{equation*}
(really there should be one finite set $T_R$ for each component of
$U_R$, but we will suppress this from the notation).

\cite{MW} shows how to extend this to an element of their $\beta\CMop(X)$
by choosing contractible open sets $V_j \subseteq X$ and embeddings
\begin{equation}
  \label{eq:18}
  \tag{*}
  D^2\times S^0\times Q_j\times V_j \to (M-\partial M) \times V_j
\end{equation}
and by taking coproducts they get an element of their $\beta\CMop(X)$
which restricts to the given germ.  To finish the proof that our
$\beta\CMtop(X)$ is contractible we have to promote \eqref{eq:18} to an
object of our $\beta\CMtop(V_j)$.  This can be done by the next lemma.

\begin{lem}
  Let $M$ be a surface and $TM \to U_2$ a bundle map.  Let $X$ be
  contractible and let $V = \R^2\times\R\times T\times X$ be the
  trivial vectorbundle over $T\times X$.  Then for any embedding
  \begin{equation*}
    e\colon  D^2 \times S^0 \times T\times X \to (M-\partial M) \times X
  \end{equation*}
  over $X$ there exists a bundle map $T^\omega V_{|\Sad(V)} \to U_3$ and
  a diffeomorphism $h\colon  M_+(0) \to D^2\times S^0\times T \times X$ such
  that the fibrewise differential of $e\circ h$ is over $U_2$.
\end{lem}
\begin{proof}
  First let $h$ be the inverse of the standard diffeomorphism given
  by~\eqref{eq:15}.  The requirement that $D(e\circ h)$ is over $U_2$
  defines a unique bundle map $T^\omega V_{|M_+(0)} \to U$, or
  equivalently a map $M_+(0) \to \Bun(\R^2\times\R,U_3)$.  After
  possibly composing $h$ with an orientation preserving diffeomorphism
  of $D^2$ we can extend this to $M_+(0) \cup (\{0\} \times D^1\times
  T \times X)$.  Now the inclusion
  \begin{equation*}
    M_+(0) \cup (\{0\}\times D^1\times T\times X) \to \Sad(V)
  \end{equation*}
  is a trivial cofibration so we can extend to all of $\Sad(V)$.
\end{proof}

\subsection{Chapter 7: Stabilisation}

This is almost as in \cite{MW}.  Start by choosing an element $z\in
\W_\emptyset(*)$ of genus 2.  This is a torus with two boundary
components and with a spin structure.  As already explained in
paragraph~\ref{sec:adapting-mw}, there are two essentially different
choices of such tori, but which one we pick is not important for
stabilisation.

As in \cite{MW} we get a fibration sequence
\begin{equation*}
  |z^{-1}h\V| \to \hocolim z^{-1} |\WT| \to \hocolim z^{-1} | \WlocT|
\end{equation*}
where $|z^{-1} h\V| \simeq \Omega^\infty\Th(-U)$ and where
\begin{equation*}
  \hofib(z^{-1}|\WT| \to z^{-1}|\WlocT|) \simeq \Z\times
  \MM^\theta(F_{\infty,2+2|T|}).
\end{equation*}
When the spaces $\MM^\theta(F_{\infty,2+2|T|})$ satisfies Harer
stability, i.e.\ if any morphism $S \to T$ in $\KK$ induces homology
equivalences $\MM^\theta(F_{\infty,2+2|T|}) \to
\MM^\theta(F_{\infty,2+2|S|})$, then the proof in \cite{MW} goes through
and proves that
\begin{equation*}
  \Z\times \MM^\infty(F_{\infty,2}) \to \Omega^\infty\Th(-U)
\end{equation*}
is a homology equivalence.

And we know from \cite{Harer} and \cite{Bauer} that for $\theta =
\theta_{\Spin}$ this Harer stability indeed does hold.

\end{document}